\documentclass[12pt]{article}

\usepackage{amsfonts}

\newtheorem{thm}{Theorem}[section]
\newtheorem{lemma}[thm]{Lemma}
\newtheorem{cor}[thm]{Corollary}

\newtheorem{prop}[thm]{Proposition}
\newtheorem{conjecture}{Conjecture}

\newenvironment{remark}{\par\medskip\noindent{\bf Remark.\ }}{\par\smallskip}
\newcommand{\proof
}{\par\medskip\noindent {\bf Proof.\ \ }}

\newcommand{\be}{\begin{equation}}
\newcommand{\ee}{\end{equation}}
\newcommand{\openbox}{\leavevmode
  \hbox to8pt{\hfil\vrule\vbox to6pt{\hrule width6pt\vfil\hrule}\vrule}}

\newcommand{\ve}[1]{\mathbf{#1}}
\newcommand{\qed}{\hbox to5pt{ } \hfill \openbox\bigskip\medskip}

\newcommand{\Fp}{\mathbb F _p}

\newcommand{\val}{\mbox{\rm val}}

\newcommand{\cF}{\mbox{$\cal F$}}

\newcommand{\N}{\mathbb N}
\newcommand{\Z}{\mathbb Z}
\newcommand{\Q}{\mathbb Q}

\newcommand{\F}{\mathbb F}

\title{New upper bounds for the size of set systems with restricted intersections modulo prime powers}
\author{G\'abor Heged\"{u}s
\\{\normalsize  \'Obuda University}
\\{\normalsize B\'ecsi \'ut 96/B, Budapest, Hungary, H-1032}
\\{\normalsize hegedus.gabor@uni-obuda.hu}
}

\begin{document}

\maketitle

\begin{abstract}
Let $q=p^\alpha$  be a fixed prime power, $k\geq 2$ be an integer. 
We give a new upper bound for the size of $k$-wise $q$-modular  $L$-avoiding  $L$-intersecting set systems, where $L$ is any proper subset of   $\{0, \ldots , q-1\}$. 
 
Our proof is based on  the linear algebra bound method and basic number theory.
\end{abstract}

\medskip

\section{Introduction}

First we introduce some notation. Let $m$ be a positive integer. Let $L\subseteq \Z$ be a finite subset. We say that $r\in L \pmod m$, if  $r\equiv \ell \pmod m$ for some $\ell\in L$. We say that   $r\notin L \pmod m$, if $r\not\equiv \ell \pmod m$ for every $\ell\in L$. 

Throughout this paper 
$n$ will be a positive integer and $[n]$ stands for the set $\{1,2,
\ldots, n\}$. The family of all subsets of $[n]$ is denoted by $2^{[n]}$. 

Let  $\F$ be a fixed field and  let 
$\F[x_1, \ldots, x_n]=\F[\ve x]$ denote  the
ring of polynomials in the variables $x_1, \ldots, x_n$ over $\F$.
For a subset $F \subseteq [n]$ we write
$\ve x_F = \prod_{j \in F} x_j$. In particular, $\ve x_{\emptyset}= 1$.

Throughout the paper  $p$ stands for a prime and 
 $\alpha> 0$ will be a positive integer. Let $q=p^\alpha$ denote a fixed prime power.

Let  $L \subseteq \{0, \ldots , q-1\}$ be a proper subset. Then the set system  $\cF$ is  $L$-avoiding mod $q$, if  $|E|\not\in L \pmod q$ for each $F\in \cF$.  The set system  $\cF$ is  $L$-intersecting mod $q$  if $|E\cap F|\in L \pmod q$ for all $E,F\in \cF, E\neq F$.

Let $m(n,s,q)$ denote the smallest integer such that for any set  $L \subseteq \{0, \ldots , q-1\}$ of size $|L|=s$, if $\cF$ is a set system of $[n]$ and $\cF$ is  $L$-avoiding mod $q$ and $L$-intersecting mod $q$  then $|\cF|\leq m(n,s,q)$. 

Deza, Frankl  and  Singhi proved an upper bound for  the size of $p$-modular $L$-avoiding  $L$-intersecting set systems of $[n]$ with $|L|=s$, in the case of $p$ prime.
\begin{thm} \label{ABS}
Let $p$ be a prime. Then  $m(n,s,p)\leq \sum_{j=0}^s {n\choose j}$.
\end{thm}

Babai et al. gave in \cite{BFKS} a polynomial upper bound, that is, of the form $O(n^{c(s)})$ for some function $c(s)$ for the size of $q$-modular $L$-avoiding  $L$-intersecting set systems of $[n]$ with $|L|=s$, when $q$ is prime power.

\begin{thm} \label{BFSK}
Let $p$ be a prime
and $\alpha> 0$ be a positive integer. Let $q=p^\alpha$ denote a prime power. Define $D:=\min (\lfloor (1+ \frac{s-1}{\alpha})^{\alpha}\rfloor, 2^{s-1})$. Then  $m(n,s,q)\leq  \sum_{j=0}^{D} {n\choose j}$.
\end{thm}

This polynomial upper bound in Theorem \ref{BFSK} fails to extend to the $m$-modular version, where $m$ is not a prime power. 
Namely Grolmusz proved in \cite{G}  that for each $m$ with at least two distinct prime divisors, there exists an   $m$-modular  $[m-1]$-avoiding  $[m-1]$-intersecting set system $\cF\subseteq 2^{[n]}$  with super-polynomial size.

Our main result is  a new  polynomial upper bound for the size of $q$-modular $L$-avoiding  $L$-intersecting set systems of $[n]$, when $q$ is prime power.

Xu and Yip proved the  following weaker bound for the size of $q$-modular $L$-avoiding  $L$-intersecting set systems of $[n]$ in \cite{XY}.
\begin{thm} \label{XY}
Let $L\subseteq \{0,1,
\ldots, q-1\}$  be a subset with size $s$. 
 If $\cF \subseteq  2^{[n]}$ is a $q$-modular $L$-avoiding  $L$-intersecting set system of $[n]$, then
$$
|\cF|\leq (q-s)\sum_{j=0}^{q-1} {n-1\choose j}.
$$
\end{thm}

Let $L\subseteq \{0,1,
\ldots, q-1\}$  be a subset with size $s>0$.  We say that a set system $\cF \subseteq  2^{[n]}$ is a $q$-modular {\em $L$-differencing Sperner}, if $|E\setminus F| \in L \pmod q$ for all $E,F\in \cF, E\neq F$. 

Xu and Yip gave the following  upper bound for the size of $q$-modular $L$-differencing Sperner set systems, when $L$ is an interval.
\begin{thm} \label{XY2}
Let $p$ be a prime
and  let $q=p^\alpha$, $\alpha\geq 1$ be a prime power.
 Let $L=\{b-s+1, \ldots , b\}$ such that $s\leq b<q$.  Assume that ${b\choose s}\not \equiv 0 \pmod p$. If $\cF \subseteq  2^{[n]}$ is a $q$-modular $L$-differencing Sperner set system, then
$$
|\cF|\leq \sum_{j=0}^s {n-1\choose j}.
$$
\end{thm}

We state our main results. 
 
First we apply Theorem \ref{XY2} to give an other new upper bound for size of $q$-modular  $L$-avoiding  $L$-intersecting set system in the special case  $L=[s]$.

\begin{thm} \label{main2}
Let $p$ be a prime
and  let $q=p^\alpha$, $\alpha\geq 1$ be a prime power. Let $1\leq s<q$ be an integer
and  let $L=\{1, \ldots , s\}$.  Let $\cF \subseteq  2^{[n]}$ be a $q$-modular $L$-intersecting set system such that $|F|\equiv 0 \pmod q$ for each $F\in \cF$. Then 
$$
|\cF|\leq \sum_{j=0}^s {n-1\choose j}.
$$
\end{thm}

Next  we give a new upper bound for size of $k$-wise $q$-modular  $L$-avoiding  $L$-intersecting set systems.

\begin{thm} \label{main3}
Let $p$ be a prime
and  let $q=p^\alpha$, $\alpha> 0$ be a prime power. 
Let $k\geq 2$ be an integer.
Let $L\subseteq \{0,1,
\ldots, q-1\}$  be a proper subset with size $s>0$.  Let $\cF \subseteq  2^{[n]}$ be a set system such that $|F|\notin L \pmod q$ for every $F\in \cF$ and $|F_1\cap \ldots \cap F_k| \in L \pmod q$ for any collection of $k$ distinct sets from $ \cF$. Then 
$$
|\cF|\leq (k-1)\sum_{j=0}^{q-1} {n\choose j}.                   
$$
\end{thm}
\begin{thm} \label{main4}
Let $p$ be a prime
and  let $q=p^\alpha$, $\alpha> 0$ be a prime power. 
Let $k\geq 2$ be an integer.
Let $L\subseteq \{0,1,
\ldots, q-1\}$  be a proper subset with size $s>0$.  Let $\cF \subseteq  2^{[n]}$ be a set system such that $|F|\notin L \pmod q$ for every $F\in \cF$ and $|F_1\cap \ldots \cap F_k| \in L \pmod q$ for any collection of $k$ distinct sets from $\cF$.
Suppose that there exists a $t\leq q-1$ integer such that  $|F|\in \{q-t, \ldots ,q-1\} \pmod q$ for each $F\in \cF$, then
$$
|\cF|\leq (k-1)\sum_{j=q-t}^{q-1} {n\choose j}.                   
$$
\end{thm}

\begin{cor} \label{main}
Let $p$ be a prime
and  let $q=p^\alpha$, $\alpha\geq 1$ be a prime power. Then 
$m(n,s,q)\leq  \sum_{j=0}^{q-1} {n\choose j}$.
\end{cor}

\begin{remark}
It is easy to verify that the upper bound of Theorem \ref{main3} is asymptotically tight, i.e. when $q$ is a fixed prime power and $n$ tends to infinity. It comes out that there exists a construction of a set system $\cF$ on $[n]$ with size 
$$
(1+o(1))(k-1){n\choose q-1}
$$ 
such that $|F|\equiv -1 \pmod q$ for every $F\in\cF$ and $|F_1\cap \ldots \cap F_k|\in L \pmod q$ for every collection of $k$ distinct sets $F_1, \ldots ,F_k\in \cF$, where $L=\{0,1,\ldots, q-2\}$. Namely we can apply the construction appearing in \cite{GS} Section 4, where we need to replace the prime $p$ with the prime power $q$ in this  construction.
\end{remark}

In Section 2 we collected all the preliminaries which we used in our proofs. In Section 3 we present the proofs of our main results.

\section{Preliminaries}

\subsection{Number theory}

The following simple fact was proved in Proposition $5.31$ of  \cite{BF}.

\begin{prop} \label{binomka} Let $q=p^{\alpha}$, $p$ a prime, and $\alpha\geq 1$.
For any integer $r$, the binomial coefficient ${r-1 \choose q-1}$ is divisible 
by $p$ iff $r$ is not divisible by $q$.  $\Box$ \end{prop}

The following result was proved in Lemma 3.3 of \cite{HR}. The congruence follows easily from the Vandermonde identity (see e.g. \cite{GKP}, pp. 169-170).
\begin{lemma} \label{binom}
Let $q=p^{\alpha}>1$ a prime power. Let $x,\ j$ be integers, $0\leq j <q$.
Then 
$$
{x+q \choose j} \equiv {x \choose j} \pmod p.
$$
\end{lemma} 

\subsection{Multilinear polynomials}

The proofs of  our main result are based on  the linear algebra bound method, some basic number theory  and the Triangular Criterion (see e.g. \cite{BF} Proposition 2.5).  We recall here shortly for the reader's convenience this basic principle.

\begin{prop} \label{tri} (Triangular Criterion)
Let $\F$ denote an arbitrary field. For $i=1,\ldots m$ let $f_i:\Omega \to \F$ be functions and $\ve v_i\in \Omega$ elements such that $f_i(\ve v_j)\neq 0$ if $i=j$ and $f_i(\ve v_j)=0$ if $i<j$. Then $f_1,\ldots ,f_m$ are linearly independent functions of the vector space $\F^{\Omega}$. 
\end{prop}

Recall that a polynomial is said to be {\em multilinear}, if it has degree at most $1$ in each variable. Let $f$ be a polynomial in $\F[x_1, \ldots, x_n]$ of degree at most $s$. Then there exists a unique  multilinear polynomial $\overline{f}$ of degree at most $s$ such that 
$$
f(\ve v)=\overline{f}(\ve v)
$$
for each $\ve v\in \{0,1\}^n$. This polynomial $\overline{f}$  is the  {\em multilinearization} of the polynomial $f$.

Finally we need for the following simple fact, which  was proved in Lemma 3.5 of \cite{HR}.

\begin{lemma} \label{egesz}
Let $f\in \Q[x_1,\ldots,x_n]$ be a polynomial such that $f(\ve v) \in \Z$ for 
every
$\ve v \in \{0,1\}^n$. 
Let $\overline{f}$ denote the multilinearization of $f$. Then $\overline{f}\in \Z[x_1,\ldots,x_n].$ 
\end{lemma}

\subsection{Separating polynomials}

We introduce first a notation. Let $I(d)\subseteq \Q[x]$ denote the set of all polynomials $p(x)$ of degree at most $d$ which take integer values $p(s)$ for every $s\in \N$.  Let $I:=\cup_{d\in \N} I(d)$, i.e. $I$ denotes the set of all polynomials $h(x)$ which take integer values $h(s)$ for every $s\in \N$. 

It can be shown easily that every  polynomial $p(x)\in I(d)$ can be written as a unique linear combination with integer coefficients of the polynomials ${x\choose 0}, \ldots , {x\choose d}$ (see \cite{BF} Ex. 7.3.3).

Now  we introduce some basic $p$-adic terminology.
Let $t$ be an integer. The $p$-adic valuation $\val(t)$ of $t$ is the exponent $j$ such that $p^j$ divides $t$, but $p^{j+1}$ does not. For example $\val(p)=1$ and $\val(1)=0$. Clearly $\val(uv)=\val(u)+\val(v)$ and $\val(u+v)\geq \min(\val(u),\val(v))$, if $u\ne 0$ and $v\ne 0$. It is easy to see that 
if $\val(u)\neq \val(v)$, then $\val(u+v)= \min(\val(u),\val(v))$.

We say that a  polynomial $h\in I$ {\em separates} a set $A\subseteq \Z$ from a set $B\subseteq \Z$ if 
$$
\max\{\val(h(s)):~ s\in A\}< \min\{\val(h(s)):~ s\in B\}.
$$

Babai et al. proved the following remarkable result in \cite{BFKS}.
\begin{lemma} \label{seppol}
Let $q=p^{\alpha}$, where $p$ is a prime and $\alpha\geq 1$.
Let $D:=\min (\lfloor (1+ \frac{s-1}{\alpha})^{\alpha}\rfloor, 2^{s-1})$. Let $L\subseteq \{0,1,
\ldots, q-1\}$  be a proper subset with size $s>0$. Then for any $\nu \in  \{0,1,
\ldots, q-1\}\setminus L$, there exists a polynomial $h\in \Z[x]$ of degree at most $D$ which separates $\nu +q\Z$ from $L+q\Z$. 
\end{lemma} 

\begin{remark}
The investigation of separating polynomials was motivated by the following Lemma (see  \cite{BFKS} Lemma 3.1).

\begin{lemma} \label{val_upper_b}
Let $p$ be a prime
and  let $q=p^\alpha$, $\alpha> 0$ be a prime power. 
Let $L\subseteq \{0,1,
\ldots, q-1\}$  be a proper subset.  Let $\cF \subseteq  2^{[n]}$ be a  $q$-modular  $L$-avoiding  $L$-intersecting  set system. 
Assume that for any $\ell\notin L \pmod q$ there exists a degree-$d$ uni-variate separating polynomial $h_{\ell}$ separating $\ell$ from $L+q\Z$. Then 
$$
|\cF|\leq \sum_{j=0}^d {n\choose j}.
$$
\end{lemma} 
\end{remark}
%

\section{Proofs}


{\bf Proof of Theorem \ref{main2}:}

 Let $\cF \subseteq  2^{[n]}$ be a $q$-modular $L$-intersecting set system such that $|F|\equiv 0 \pmod q$ for each $F\in \cF$. Define $L':=\{q-1, \ldots , q-s\}$. It is easy to verify that $\cF$ is a $q$-modular $L'$-differencing Sperner set system. Since $p\not| {q-1 \choose s}$, hence it follows from Theorem \ref{XY2} that 
$$
|\cF|\leq \sum_{j=0}^s {n-1\choose j}.
$$
\qed

The proof of Theorem \ref{main3} is based on our following construction of $q$ modular separating polynomials.

\begin{prop} \label{interpol}
Let $q=p^{\alpha}>1$ a prime power. Let $L \subseteq \{0, \ldots , q-1\}$ be a fixed proper subset. Then there exists an uni-variate  polynomial $F_L(x)\in \Q[x]$ such that \\
\begin{itemize}
\item[(1)] $F_L(t)\in \Z$ for each $t\in \Z$;

\item[(2)] $F_L(t)$ separates $(\{0, \ldots , q-1\}\setminus L)+q\Z$  from  $L+q\Z$;

\item[(3)] $F_L(m)\equiv 1 \pmod p$ for each  $m\notin L \pmod q$ and $F_L(s)\equiv 0 \pmod p$ for each $s\in L \pmod q$  and 

\item[(4)] $\deg(F_L)\leq q-1$.
\end{itemize}
\end{prop} 
\proof 

Let $a\in \{0, \ldots , q-1\}$ be fixed. Consider the polynomial 
$$
L_a(x):=(-1)^{q-1} {x-a-1\choose q-1}\in \Q[x].
$$
It follows from Proposition \ref{binomka} that  $L_a(a)\equiv 1 \pmod p$ and $L_a(b)\equiv 0 \pmod p$ for each $b\in \{0, \ldots , q-1\}\setminus \{a\}$. 

Define
$$
F_L(x):=\sum_{a\in \{0, \ldots , q-1\}\setminus L} L_a(x)\in \Q[x].
$$
It follows from the definition of $F_L$ that $F_L(t)\in \Z$ for each $t\in \Z$.

Clearly $F_L(m)\equiv 1 \pmod p$ for each  $m\in \{0, \ldots , q-1\}\setminus L$ and $F_L(m)\equiv 0 \pmod p$ for each $m\in L$. It follows from Lemma \ref{binom} that
$$
F_L(t+mq)\equiv F_L(t) \pmod p
$$
for each $m.t\in \Z$. Consequently $F_L(m)\equiv 1 \pmod p$ for each  $m\notin L \pmod q$ and $F_L(s)\equiv 0 \pmod p$ for each $s\in L \pmod q$. This means that $F_L(t)$ separates  $(\{0, \ldots , q-1\}\setminus L)+q\Z$ from $L+q\Z$. 

We get $\deg(F_L)\leq q-1$ from the definition of $F_L$. \qed


{\bf Proof of Theorem \ref{main3}:}
                     
Let  $L \subseteq \{0, \ldots , q-1\}$ be a proper subset and let $\cF$ be a family of sets satisfying the  conditions of Theorem \ref{main3}.  We will repeat the following natural procedure until $\cF$ becomes empty. 
                    
Suppose that $\cF\ne \emptyset$ at round $i$. Then we can choose a maximal collection of sets $F_1, \ldots , F_d$ from the family $\cF$ such that $|\cap_{j=1}^e F_j| \notin L\pmod q$ for each $1\leq e\leq d$, but for any additional new set $F'\in \cF$ we get that $|\cap_{j=1}^d F_j\cap F'| \in L\pmod q$. Then define $X_i:=F_1$ and $Y_i:=\cap_{j=1}^d F_j$. Finally we clear all the sets $F_1, \ldots , F_d$ from the family $\cF$. Clearly after this process we get at least $m\geq \frac{|\cF|}{k-1}$ pairs of subsets ($X_i$, $Y_i$). Next we derive an upper bound for $m$.

It follows from the definition that $|X_i\cap Y_i|=|Y_i| \notin L \pmod q$, but $|X_t\cap Y_i| \in L \pmod q$ for each $t>i$. Let  $\ve x_i$ and $\ve y_i$ denote the characteristic vectors of the sets  $X_i$, $Y_i$, respectively. 

We denote by $\ve x\cdot \ve y:=\sum_i x_iy_i$ the standard scalar product for any two vectors $\ve x, \ve y\in \Q^n$. Obviously $\ve x_t \cdot \ve y_i=|X_t \cap Y_i|$ for any $1\leq t,i\leq m$. 

Define the polynomial in $n$ variables
$$
F_i(\ve x):= F_L(\ve x \cdot \ve y_i)\in \Q[\ve x]
$$
for each $1\leq i\leq m$. Recall that we defined the polynomial  $F_L$ in Proposition \ref{interpol}. 

It follows from Lemma  \ref{egesz} that the multilinear polynomial  $\overline{F_i}\in \Z[\ve x]$.

Then define $G_i\in {\Fp}[\ve x]$ as the modulo $p$ reduction of the multilinear polynomial  $\overline{F_i}\in \Z[\ve x]$. Clearly $\deg(G_i)\leq \deg(F_i)\leq q-1$ for each $1\leq i\leq m$.

It is easy to verify that
$$
G_i(\ve x_i)=F_i(\ve x_i)=F_L(|X_i \cap Y_i|)=F_L(|Y_i|) \pmod p,
$$
hence $G_i(\ve x_i)\not\equiv 0 \pmod p$ for each $i$, because $|Y_i|=|\cap_{j=1}^d F_j|\notin L \pmod q$.

On the other hand 
$$
G_i(\ve x_t)=F_i(\ve x_t)=F_L(|X_t \cap Y_i|) \pmod p,
$$
which means that $G_i(\ve x_t)\equiv 0 \pmod p$ for each $1\leq i<t\leq m$, because $|X_t\cap Y_i| \in L \pmod q$ for each $t>i$. 

It follows from Proposition \ref{tri}  (Triangular Criterion) that the    polynomials  $\{G_i:~ 1\leq i\leq m\}$ are  linearly independent.  Clearly each $G_i$ belongs to the space of multilinear   polynomials of degree at most $q-1$. Since the dimension of this vector space is $\sum_{j=0}^{q-1} {n\choose j}$, hence we get the desired upper bound on  $m$, which implies the bound on $|\cF|$. \qed   

{\bf Proof of Theorem \ref{main4}:}
We prove here Theorem \ref{main4}  combining  the linear algebra bound method with the nice technique by Blokhuis (see e.g. \cite{B}).

First we define a new set of multilinear  polynomials. Consider the polynomials
$$
k_I(\ve x):= x_I\cdot \prod_{\ell=1}^t \Big( \sum_{j=1}^n x_j-(q-\ell) \Big)\in \Fp[\ve x]
$$
for each $I\subseteq [n]$. Let $h_I(\ve x):=\overline{k_I}(\ve x)$ for $I\subseteq [n]$. It is not difficult to verify that polynomials  $\{h_I:~ |I|\leq q-1-t\}$ are  linearly independent over $\Fp$. For the proof of this fact we may
use the condition that $|F|\in \{q-t, \ldots ,q-1\} \pmod q$ for each $F\in \cF$.

Finally we show that the set of polynomials  
$$
\{G_i:~ 1\leq i\leq m\}\cup \{h_I:~ |I|\leq q-1-t\}
$$ 
is linearly independent over $\Fp$. 
Assume the contrary, that the  polynomials $\{G_i:~ 1\leq i\leq m\}\cup \{h_I:~ |I|\leq q-1-t\}$  are not linearly independent.

Consider a  non-trivial linear combination 
\begin{equation}  \label{lincomb} 
\sum_{i=1}^m \alpha_i G_i(\ve x) + \sum_{I:|I|\leq q-t-1}  \beta_I h_I(\ve x)=0
\end{equation}
for some $\alpha_i,\beta_I \in \Fp$. 

Since the polynomials  $\{h_I:~ |I|\leq q-1-t\}$ are  linearly independent,  
the first sum $\sum_{i=1}^m \alpha_i G_i(\ve x)$ itself  is a non-trivial linear combination. 

Let $i_0$ be largest index such that $\alpha_{i_0}\ne 0$.  Then by substituting $\ve x:=\ve x_{i_0}$ all terms in the  second sum $\sum_{I:|I|\leq q-t-1}  \beta_I h_I(\ve x)$ vanish, because $|F|\in \{q-t, \ldots ,q-1\} \pmod q$  for each $F\in \cF$ and hence $h_I(\ve x_{i_0})=0 \pmod p$ for each $I$. Consequently  all terms of the first sum but the one with the index $i_0$ vanish. But $G_{i_0}(\ve x_{i_0})\not \equiv 0 \pmod p$, which implies that  $\alpha_{i_0}= 0$, a contradiction.

Let $V$ denote the vector space of multilinear polynomials in $n$ variables of degree at most $q-1$ over $\Fp$. Obviously $\dim_{\Fp} V=\sum_{i=0}^{q-1} {n \choose i}$. 

We have found $m+\sum_{i=0}^{q-t-1} {n \choose i}$  linearly independent  polynomials in $V$. Consequently 
$$
\frac{|\cF|}{k-1}\leq m\leq \dim_{\Fp} V - \sum_{i=0}^{q-t-1} {n \choose i}=\sum_{i=q-t}^{q-1} {n \choose i}.
$$
\qed

\section{Concluding remarks}
We propose the following conjecture as a strengthening of Theorem \ref{main3}.

\begin{conjecture} \label{Hconj}
Let $p$ be a prime
and  let $q=p^\alpha$, $\alpha\geq 1$ be a prime power. Let $k\geq 2$ be an integer.
Let $L\subseteq \{0,1,
\ldots, q-1\}$  be a proper subset with size $s>0$.  Let $\cF \subseteq  2^{[n]}$ be a set system such that $|F| \notin L  \pmod q$ for every $F\in \cF$ and $|F_1\cap \ldots \cap F_k| \in L \pmod q$ for any collection of $k$ distinct sets from $ \cF$. 
Suppose that there exist $t\leq q-1$  integers $k_1, \ldots ,k_t\in  \{0,1,
\ldots, q-1\}$ such that   $|F|\pmod q\in \{k_1, \ldots ,k_t\}$ for each $F\in \cF$. Then
$$
|\cF|\leq (k-1)\sum_{j=q-t}^{q-1} {n\choose j}.                   
$$
\end{conjecture}

It is not very difficult to verify that combining our proof of Theorem \ref{main3} and  Theorem \ref{main4} with Lemma \ref{seppol} gives the following result. The proof of this statement will be left to the reader.

\begin{thm} \label{main5}
Let $q=p^{\alpha}$, where $p$ is a prime and $\alpha\geq 1$.
Let $k\geq 2$ be an integer.
Let $L\subseteq \{0,1,
\ldots, q-1\}$  be a subset with size $s>0$. Let $D:=\min (\lfloor (1+ \frac{s-1}{\alpha})^{\alpha}\rfloor, 2^{s-1})$.  Let $\cF \subseteq  2^{[n]}$ be a set system such that $|F| \notin L  \pmod q$ for every $F\in \cF$ and $|F_1\cap \ldots \cap F_k| \in L \pmod q$ for any collection of $k$ distinct sets from $ \cF$.  Then 
$$
|\cF|\leq (k-1)\sum_{j=0}^{D} {n\choose j}.                   
$$
If in addition there exists a  $t\leq D$ integer  such that  $|F|\in \{D-t+1, \ldots ,D\}\pmod q$ for each $F\in \cF$, then
$$
|\cF|\leq (k-1)\sum_{j=D-t+1}^{D} {n\choose j}.                   
$$
\end{thm}

{\bf Acknowledgement.}           
I am indebted to  Lajos R\'onyai for his useful remarks.

\end{document}